\documentclass[10pt,leqno]{article}
\usepackage{amssymb,amsbsy,amsmath,amsfonts,amssymb,amscd, mathrsfs}
\usepackage{graphicx} 

\usepackage[english]{babel}
\usepackage[T1]{fontenc}
\usepackage{indentfirst}

\makeatletter
\@addtoreset{equation}{section}
\makeatother

\newtheorem{statement}{}[section]
\newtheorem{theorem}[statement]{Theorem}
\newtheorem{lemma}[statement]{Lemma}

\newtheorem{definition}[statement]{Definition}

\newtheorem{claim}[statement]{Claim}

\newcommand\C{\mathbb C}

\newcommand\R{\mathbb R}
\newcommand\T{\mathbb T}
\newcommand\D{\mathbb D}

\renewcommand\H{\mathbb H}
\newcommand\e{{\rm e}}

\newcommand\eps{\varepsilon}
\newcommand\ind{{\rm 1\kern-.30em I}}
\newcommand\qed{\hfill $\square$}
\renewcommand \Re{{\mathfrak R}{\rm e}\,}
\renewcommand \Im{{\mathfrak I}{\rm m}\,}

\let\phi=\varphi

\title{\bf Estimates for approximation numbers of some classes of composition operators on the Hardy space}
\author{\it Daniel Li, Herv\'e Queff\'elec, Luis Rodr{\'\i}guez-Piazza\footnote{Supported by a Spanish research project MTM 2009-08934.}}

\date{\footnotesize \today}

\begin{document}

\maketitle

\noindent{\bf Abstract.} \emph{We give estimates for the approximation numbers of composition operators on $H^2$, in terms of some 
modulus of continuity. For symbols whose image is contained in a polygon, we get that these approximation numbers are dominated by $\e^{- c \sqrt n}$. 
When the symbol is continuous on the closed unit disk and has a domain touching the boundary non-tangentially at a finite number of points, with a good behavior 
at the boundary around those points, we can improve this upper estimate. A lower estimate is given when this symbol has a good radial behavior at some 
point. As an application we get that, for the cusp map, the approximation numbers are equivalent, up to constants, to $\e^{- c \, n / \log n }$, very near to the 
minimal value $\e^{- c \, n}$. We also see the limitations of our methods. To finish, we improve a result of O. El-Fallah, K. Kellay, M. Shabankhah and 
H. Youssfi, in showing that for every compact set $K$ of the unit circle $\T$ with Lebesgue measure $0$, there exists a compact composition operator 
$C_\phi \colon H^2 \to H^2$, which is in all Schatten classes, and such that $\phi = 1$ on $K$ and $|\phi | < 1$ outside $K$. }
\medskip

\noindent{\bf Mathematics Subject Classification 2010.} Primary: 47B06 -- Secondary: 30J10 ; 47B33
\medskip

\noindent{\bf Key-words.} approximation numbers; Blaschke product; composition operator; cusp map; Hardy space; modulus of continuity; Schatten classes   


\section{Introduction and notation} 

If the approximation numbers of some classes of operators on Hilbert spaces are well understood (for example, those of Hankel operators: see \cite{MPT}), 
it is not the case of those of composition operators. Though their behavior remains mysterious, some recent results are obtained  in \cite{LIQUEROD} and 
\cite{lens} for approximation numbers of composition operators on the Hardy space $H^2$. In \cite{LIQUEROD}, it is proved that one always has 
$a_n (C_\phi) \gtrsim \e^{- c \, n}$ for some $c > 0$ (\cite{LIQUEROD}, Theorem~3.1) and that this speed of decay can only be got when the symbol $\phi$ 
maps the unit disk $\D$ into a disk centered at $0$ of radius strictly less than $1$, i.e. $\|\phi\|_\infty < 1$ (\cite{LIQUEROD}, Theorem~3.4). \par 

In this paper, we give estimates which are somewhat general, in terms of some modulus of continuity. In Section~\ref{disk algebra}, we obtain an upper estimate 
when the symbol $\phi$ is continuous on the closed unit disk and has an image touching non-tangentially the unit circle at a finite number of points, with a good 
behavior on the boundary around this point. As an application, we show that for symbols $\phi$ whose image is contained in a polygon 
$a_n (C_\phi) \leq a\, \e^{- b \sqrt n}$, for some constants $a, b > 0$; this has to be compared with \cite{lens}, Proposition~2.7, where it is shown that if 
$\phi$ is a univalent symbol such that $\phi (\D)$ contains an angular sector centered on the unit circle and with opening $\theta \pi$, $0 < \theta < 1$, then 
$a_n (C_\varphi) \geq a \, \e^{- b \sqrt n}$, for some (other) positive constants $a$ and $b$, depending only on 
$\theta$. In Section~\ref{lower}, we obtain a lower bound when $\phi$ has a good radial behavior at the contact point. Both proofs use 
Blaschke products. This allows to recover the estimation $a_n (C_{\lambda_\theta}) \approx \e^{ - c \, \sqrt n}$ obtained in \cite{LIQUEROD}, 
Proposition~6.3, and \cite{lens}, Theorem~2.1 for the lens map $\lambda_\theta$. In Section~\ref{section cusp}, we give another example, the cusp map, for 
which $a_n (C_\phi) \approx \e^{- c \, n / \log n}$, very near the minimum value $\e^{- c \, n}$. We end that section by considering a one-parameter class of 
symbols, first studied by J. Shapiro and P. D. Taylor \cite{SHTA} and seeing the limitations of our methods. In Section~\ref{application}, we improve a result of 
E.A.~Gallardo-Guti\'errez and M.J.~Gonz\'alez (previously generalized by O. El-Fallah, K. Kellay, M. Shabankhah and H. Youssfi \cite{EKSY}, Theorem~3.1). 
It is known that for every compact composition operator $C_\phi \colon H^2 \to H^2$, the set 
$E_\phi = \{ \e^{i \theta} \, ; \ |\phi^\ast ( \e^{i \theta}) | = 1 \}$ has Lebesgue measure $0$. These authors showed (\cite{GALGON}), with a rather 
difficult construction, that there exists a compact composition operator $C_\phi \colon H^2 \to H^2$ such that the Hausdorff dimension of $E_\phi$ is equal to 
$1$ (and in \cite {EKSY}, it is shown that for any negligible compact set $K$, there is a Hilbert-Schmidt operator $C_\phi$ such that $E_\phi = K$). We 
improve this result in showing that for every compact set $K$ of the unit circle $\T$ with Lebesgue measure $0$, there exists a compact composition operator 
$C_\phi \colon H^2 \to H^2$, which is even in all Schatten classes, and such that $E_\phi = K$. 
\medskip

\noindent \emph {Notation.} We denote by $\D$ the open unit disk and by $\T = \partial \D$ the unit circle; $m$ is the normalized Lebesgue measure on $\T$: 
$dm (t) = dt / 2 \pi$. The disk algebra $A (\D)$ is the space of functions which are continuous on the closed unit disk $\overline{\D}$ and analytic in 
the open unit disk. If $H^2$ is the usual Hardy space on $\D$, every analytic self-map $\phi \colon \D \to \D$ (also called \emph{Schur function}) defines, 
by Littlewood's subordination principle, a bounded operator $C_\phi \colon H^2 \to H^2$ by $C_\phi (f) = f \circ \phi$, called the 
\emph{composition operator of symbol $\phi$}. 
\par\medskip

Recall that if $T \colon E \to F$ is a bounded operator between two Banach spaces, the \emph{approximation numbers $a_n (T)$} of $T$ are defined by:
\begin{displaymath} 
\qquad \quad a_n (T) = \inf \{ \| T - R\| \, ; \ {\rm rank}\, (R) < n\} \,, \qquad n = 1 ,  2 , \ldots
\end{displaymath} 
The sequence $\big( a_n (T) \big)_n$ is non-increasing and, when $F$ has the Approximation Property, $T$ is compact if and only if $a_n (T)$ tends to $0$. 
\par\smallskip

\begin{definition}
A \emph{modulus of continuity} $\omega$ is a continuous function 
\begin{displaymath} 
\omega \colon [0, A] \to \R^+, 
\end{displaymath} 
which is increasing, sub-additive, and vanishes at zero. 
\end{definition}

Some examples are: 
\begin{displaymath} 
\omega (h) = h^{\alpha},\quad  0 < \alpha \leq 1 \,; \qquad  \omega (h) = h \log \frac{1}{ h} \,;  \qquad \omega (h) = \frac{1}{ \log\frac{1}{h}} \, \cdot
\end{displaymath} 
For any modulus of continuity $\omega$, there is a concave modulus of continuity $\omega'$ such that $\omega \leq \omega' \leq 2 \omega$ (see \cite{concave} 
for example); therefore we may and shall assume that $\omega$ is concave on $[0, A]$. In that case, $\omega^{- 1}$ is convex, and 
\begin{equation} \label{defi} 
r_\omega (x) := \frac{\omega^{- 1} (x)}{x} 
\end{equation}
is non-decreasing.  
\medskip

The notation $u (t) \lesssim v (t)$ means that $u (t) \leq A v (t) $ for some constant $A > 0$ and $u (t) \approx v (t)$ means that both $u (t) \lesssim v (t)$ and 
$v (t) \lesssim u (t)$.

\section{Upper bound and boundary behavior} \label{disk algebra} 

\begin{definition} \label{def regul}
Let $\omega$ be a modulus of continuity and $\phi$ a symbol in the disk algebra $A (\D)$. Let $\xi_0 \in \partial \D \cap \phi (\overline \D)$. 
We say that the symbol $\phi$ has an \emph{$\omega$-regular behavior} at $\xi_0$ if, setting: 
\begin{equation} 
\gamma (t) = \phi (\e^{i t}) \,,
\end{equation} 
and $E_{\xi_0} = \{ t \, ; \ \gamma (t) = \xi_0 \}$, there exists $r_0 > 0$ such that:
\par\smallskip

1) for some positive constant $C > 0$, one has, for every $t_0 \in E_{\xi_0}$ and $|t - t_0| \leq r_0 $:
\begin{equation} \label{reluc} 
|\gamma (t) - \gamma (t_0) | \leq C \big(1 - | \gamma (t) | \big)  \,. 
\end{equation}
\par 

2) for some positive constant $c > 0$, one has, for for every $t_0 \in E_{\xi_0}$ and $|t - t_0| \leq r_0$: 
\begin{equation} \label{ponct}
c \, \omega (| t - t_0 |) \leq | \gamma (t) - \gamma (t_0) |  \, .  
\end{equation} 
\end{definition}
\smallskip

The first condition implies that the image of $\phi$ touches $\partial \D$ at the point $\xi_0$, and non-tangentially. The second one implies that $\phi$ does not 
stay long near $\xi_0 = \gamma (t_0)$. \par

Note that, due to \eqref{ponct}, the intervals $[t - r_0/2, t + r_0/2]$, for $t \in E_{\xi_0}$ are pairwise disjoint and therefore the set $E_{\xi_0}$ must be finite. 
\par
\medskip  

We shall make the following assumption (to avoid the Lipschitz class): 
\begin{equation} \label{assu} 
\qquad \lim_{h\to 0^+} \frac{\omega (h)}{h} = \infty \, ; \quad \text{equivalently} \quad \lim_{h \to 0^+} \frac{\omega^{- 1} (h)}{h} = 0. 
\end{equation}

Indeed, assume that $\gamma$ is $K$-Lipschitz at some point $t_0 \in [0, 2\pi]$, namely $| \phi (\e^{i t}) - \phi (\e^{i t_0}) | \leq K\, |t - t_0|$, with 
$|\phi (\e^{i t_0}) | = 1$; then 
\begin{align*} 
m (\{ t \in [0, 2\pi] \, ; \, | \phi (\e^{i t}) & - \phi (\e^{i t_0}) | \leq h\}) \\
& \geq m ( \{ t \in [0, 2\pi] \, ; \, |t - t_0| \leq h / K \}) = h / 2\pi K \,; 
\end{align*} 
hence this measure in not $o\, (h)$ and the composition operator $C_\phi$ is not compact (\cite{McC}, or \cite{Co-McC}, Theorem~3.12). \par
\medskip

In order to treat the case where the image of $\phi$ is a polygon, we need to generalize the above definition. We ask not only that $\phi$ is $\omega$-regular 
at the points $\xi_1, \ldots , \xi_p$ of contact of $\phi (\overline \D)$ with $\partial \D$, but a little bit more. \par
\begin{definition} \label{def glob-reg}
Assume that $\phi (\overline \D) \cap \partial \D = \{\xi_1, \ldots, \xi_p \}$. We say that $\phi$ is \emph{globally-regular} if there exists a modulus of 
continuity $\omega$ such that, writing $E_{\xi_j} = \{ t \, ; \ \gamma (t) = \xi_j \}$, one has, for some $r_1, \ldots, r_p > 0$
\begin{displaymath} 
\T = \bigcup_{j = 1}^p \big( E_{\xi_j} + [-r_j, r_j] \big)
\end{displaymath} 
and for some positive constants $C, c > 0$, 
\par\smallskip

1') one has, for $j = 1, \ldots, p$, every $t_j \in E_{\xi_j}$ and $|t - t_j| \leq r_j $:
\begin{equation} \label{reluc glob} 
|\gamma (t) - \gamma (t_j) | \leq C \big(1 - | \gamma (t) | \big)  \,. 
\end{equation}
\par 

2') one has, for $j = 1, \ldots, p$, every $t_j \in E_{\xi_j}$ and $|t - t_j| \leq r_j$: 
\begin{equation} \label{ponct glob}
c \, \omega (| t - t_j |) \leq | \gamma (t) - \gamma (t_j) |  \, .  
\end{equation} 
\end{definition}
\par

Let us note that condition 1') is equivalent to say that $\phi (\overline \D)$ is contained in a polygon inside $\overline \D$ whose vertices contain 
$\xi_1, \ldots, \xi_p$, and these are the only vertices in the boundary $\partial \D$. Of course, we may assume that \eqref{reluc glob} and \eqref{ponct glob} 
hold only when $t$ is in a neighborhood of $t_j$, since they will then hold for $|t - t_j| \leq r_j$, provided  we change the constants $C, c$. \par
\medskip

Before stating our theorem, let us introduce a notation. If $\phi$ is as in Definition~\ref{def glob-reg} and $\sigma, \kappa > 0$ are some constants, we set:
\begin{equation} \label{d_N}
d_N = \bigg[ \sigma \log \frac{ \kappa \, 2^{- N}}{\omega^{- 1} (\kappa \, 2^{- N})} \bigg] + 1 ,
\end{equation} 
where $[\ \ ]$ stands for the integer part. For every integer $q \geq 1$, we denote by 
\begin{equation} \label{N_q}
N = N_q \quad \text{the largest integer such that } p \, N d_N < q 
\end{equation}
($N_q = 1$ if no such $N$ exists). \par
\smallskip

We then have the following result.

\begin{theorem} \label{general}
Let $\phi$ be a symbol in $A (\D)$ whose image touches $\partial\D$ at the points $\xi_1, \ldots, \xi_p$, and nowhere else. Assume that $\phi$ is 
globally-regular. Then, there are  constants $\kappa$, $K$, $L > 0$, depending only on $\phi$, such that, using the notation \eqref{d_N} and \eqref{N_q}, 
one has, for every $q \geq 1$: 
\begin{equation} \label{surprise} 
a_q (C_\phi) \leq K \, \sqrt{\frac{\omega^{- 1} (\kappa \, 2^{- N_q})}{\kappa \, 2^{- N_q}}} \,. 
\end{equation}
\end{theorem}
\smallskip

Before proving this theorem, let us indicate two applications. In these examples, we can give an upper estimate for all approximation numbers $a_n (C_\phi)$, 
$n \geq 1$ because we can interpolate between the integers $N d_N$ and $(N + 1) \, d_{N + 1}$, which is not the case in general. \par \smallskip

1) $\omega (h) = h^{\theta}$, $0 < \theta < 1$, as this is the case for inscribed polygons (see the proof of the foregoing Theorem~\ref{polygon}; here 
$\theta = \max \{\theta_1, \ldots, \theta_p \}$, where $\theta_1 \pi, \ldots, \theta_p \pi$ are the values of the angles of the polygon),  as well as, with $p = 2$, 
for lens maps $\lambda_\theta$ (see \cite{Shap-livre}, page~27, for the definition; see also \cite{lens}). We have here $\omega^{- 1} (h) = h^{1 / \theta}$. 
Hence $d_N \approx N$, $N_q \approx \sqrt q$, and we then get from \eqref{surprise} that $a_q  (C_\phi) \leq \alpha \, 2^{- \delta N}$ for $q \gtrsim N^2$, 
with $\delta > 0$. Equivalently, for suitable constants $\alpha, \beta > 0$, 
\begin{equation} \label{majo lens}
a_n (C_\phi) \leq \alpha \, \e^{- \beta \sqrt n} \,,
\end{equation} 
which is the result obtained  in \cite{lens}, Theorem~2.1.\par
\smallskip

2) $\omega (h) = \frac{1}{( \log 1 /h )^\alpha}\,$\raise 1pt \hbox{,} $0 < \alpha \leq 1$, as this is the case, when $\alpha = 1$, for the cusp map, 
defined below in Section~\ref{section cusp} (with $p = 1$). Then, we have $\omega^{- 1} (h) = \e^{- h^{- 1 / \alpha}}$ and $d_N \approx 2^{N /\alpha}$, 
so that $N_q \approx \log q$ and $2^{N_q / \alpha} \approx q/ \log q$. Now, a simple computation gives: 
\begin{equation} \label{majo cusp}
a_n (C_\phi) \leq \alpha \, \e^{- \beta n /\log n} \,.
\end{equation} 
\bigskip \goodbreak 

Without assuming some regularity, one has the following general upper estimate.

\begin{theorem} \label{polygon} 
Let $\phi \colon \D \to \D$ be an analytic self-map whose image is contained in a polygon ${\mathbf P}$ with vertices on the unit circle. Then, there exist 
constants $\alpha, \beta > 0$, $\beta$ depending only on ${\mathbf P}$, such that:
\begin{equation}\label {two} 
a_{n} (C_\phi) \leq \alpha \, \e^{- \beta \sqrt{n}}. 
\end{equation}
\end{theorem} 

In \cite{lens}, Proposition~2.7, it is shown that if $\phi$ is a univalent symbol such that $\phi (\D)$ contains an angular sector centered on the unit circle and 
with opening $\theta \pi$, $0 < \theta < 1$, then $a_n (C_\varphi) \geq \alpha \, \e^{- \beta \sqrt n}$, for some (other) positive constants $\alpha$ and 
$\beta$, depending only on $\theta$. Note that the injectivity of the symbol is there necessary, since there exists (see the proof of Corollary~5.4 in 
\cite{LIQUEROD}), for every sequence $(\eps_n)$ of positive numbers tending to $0$, a symbol $\phi$ whose image is $\D \setminus \{0\}$, and hence 
contains polygons), which is $2$-valent, and for which $a_n (C_\phi) \lesssim \e^{- \eps_n n}$. This bound may be much smaller than $\e^{- \beta \sqrt n}$.
\par\bigskip 

\noindent{\bf Proof of Theorem~\ref{general}.} It follows the lines of that of \cite{lens}, Theorem~2.1. \par 
Recall (\cite{lens}, Lemma~2.4) that for every Blaschke product $B$ with less than $N$ zeros (each of them being counted with its multiplicity), one has: 
\begin{equation} \label{majo avec Blaschke}
\big[ a_N (C_\phi) \big]^2 \lesssim \sup_{0 < h < 1, |\xi| = 1} \frac{1}{h} \int_{S (\xi, h)} |B (z)|^2 \, dm_\phi (z) \,,
\end{equation} 
where $S (\xi, h) = \{ z \in \overline{\D} \, ; \ |z - \xi| \leq h\}$ and $m_\phi$ is the pull-back measure by $\phi$ of the normalized Lebesgue measure 
$m$ on $\T$. 
\smallskip

The proof will come from an adequate choice of a Blaschke product. \par

Fix a positive integer $N$. \par

Set, for $j = 1, \ldots, p$ and $k = 1, 2, \ldots$:
\begin{equation} 
p_{j, k} = (1 - 2^{- k}) \xi_j
\end{equation} 
and consider the Blaschke product of length $p Nd$ ($d$ being a positive integer, to be specified later) given by:
\begin{equation} 
B (z) = \prod_{j = 1}^p \prod_{k = 1}^N \bigg[ \frac{z - p_{j, k}}{1 - \overline{p_{j, k}} \, z} \bigg]^d .
\end{equation} 

Recall that we have set 
\begin{equation} 
\gamma (t) = \phi (\e^{i t}) . 
\end{equation} 

To use \eqref{majo avec Blaschke}, note that if $| \gamma (t) - \xi| \leq h$, then, for some $j = 1, \ldots, p$ and some $t_j \in E_{\xi_j}$, one has 
$|t - t_j| \leq r_j$ and, by \eqref{reluc glob}, $|\gamma (t) - \xi_j| \leq C ( 1 - |\gamma (t)|) \leq C \,|\gamma (t) - \xi| \leq C h$. Therefore, 
denoting by $L_j$ the number of elements of $E_{\xi_j}$ (which is finite by the remark following Definition~\ref{def regul}):
\begin{displaymath} 
\big[ a_N (C_\phi) \big]^2 \lesssim \sup_{0 < h < 1} \frac{1}{h} 
\sum_{j = 1}^p L_j \int_{\{| \gamma (t) - \xi_j | \leq C h \}\cap \{|t - t_j| \leq r_j\}} |B [\gamma (t)]|^2 \, \frac{dt}{2 \pi} \,,
\end{displaymath} 
and we only need to majorize the integrals:
\begin{displaymath} 
I_j (h) = \int_{\{| \gamma (t) - \xi_j | \leq C h \}\cap \{|t - t_j| \leq r_j\}} | B\big( \gamma (t) \big) |^2 \, \frac{dt}{2 \pi} \, \cdot
\end{displaymath} 
Moreover, it suffices, by interpolation, to do that with $h = h_n$, where $h_n = 2^{- n}$. \par 
\smallskip

By \eqref{ponct glob}, for $|t - t_j| \leq r_j$ and $| \gamma (t) - \xi_j| \leq C h_n$, one has 
$c \, \omega (|t - t_j|) \leq | \gamma (t) - \xi_j | \leq C h_n = C \, 2^{- n}$, which implies that 
\begin{equation} \label {gogo} 
| t  - t_j| \leq \omega^{- 1} (c^{- 1} C \, 2^{- n}) . 
\end{equation} 
Let
\begin{equation} 
s_n = \omega^{- 1} (c^{- 1} C \, 2^{- n}) .  
\end{equation} 
One has: 
\begin{displaymath} 
I_j (h_n) \leq  \int_{\{|t - t_j| \leq s_n \} \cap \{ |t - t_j| \leq r_j \}} | B \big( \gamma (t) \big) |^2 \, \frac{dt}{2 \pi} \, \cdot
\end{displaymath} 
For $n \geq N$, we simply majorize $| B \big( \gamma (t) \big) |$ by $1$ and we get: 
\begin{align*} 
\frac{1}{h_n} \, I_j (h_n) 
& \leq \frac{1}{h_n} \, \frac{2 s_n}{2 \pi}  
= \frac{c^{- 1} C}{\pi}\, \frac{1}{c^{- 1} C \, 2^{- n}} \, \omega^{- 1}(c^{- 1} C \, 2^{- n})  \\ 
& \leq   \frac{c^{- 1} C}{\pi}\, \frac{\omega^{- 1} (c^{- 1} C \, 2^{- N})}{c^{- 1} C \,  2^{- N}}  \, \raise 1pt \hbox{,}
\end{align*} 
since the function $\omega^{- 1} (x) / x$ is non-decreasing.  \par \smallskip 

When $n \leq N - 1$, we write:
\begin{align*} 
I_j (h_n) 
& \leq  \int_{\{|t - t_j| \leq s_N \} \cap \{|t - t_j| \leq r_j \}}  | B \big( \gamma (t) \big) |^2 \, \frac{dt}{2 \pi}  \\
& \hskip 50 pt +  \int_{\{ s_N < |t - t_j| \leq s_n \} \cap \{ |t - t_j| \leq r_j\}}   
| B \big( \gamma (t) \big) |^2 \, \frac{dt}{2 \pi} \, \cdot
\end{align*} 
The first integral is estimated as above. For the second one, we claim that:
\begin{claim}
For some constant $\chi < 1$, one has, for $j = 1, \ldots, p$ and every $t_j \in E_{\xi_j}$:
\begin{equation} 
| B \big( \gamma (t) \big)| \leq \chi^d \qquad \text{when } |t - t_j| > s_N \text{ and } |t - t_j| \leq r_j.
\end{equation} 
\end{claim}  

To see that, we shall use \cite{lens}, Lemma~2.3. Let us recall that this lemma asserts that for $w, w_0 \in\D$ satisfying 
$\vert w - w_0 \vert \leq M \min (1 - \vert w \vert, 1 - \vert w_0 \vert)$ for some positive constant $M$, one has:
\begin{equation} \label{petit lemme}
\bigg| \frac{w - w_0} {1 - \overline{w_0}\, w} \bigg| \leq \frac{M}{\sqrt{M^2 +1}}  \, \cdot
\end{equation} 

Let $t$ such that $|t - t_j| \leq r_j$ and $|t - t_j| > s_N$. We have, on the one hand, $\omega (|t - t_j|) \geq \omega (s_N) = c^{- 1} C \, 2^{- N}$, 
and, on the other hand, since $|\gamma (t_j)| = |\xi_j| = 1$
\begin{displaymath} 
c \, \omega (|t - t_j|) \leq |\gamma (t) - \gamma (t_j)| \leq C (1 - |\gamma (t)|) \,;
\end{displaymath} 
hence $1 - |\gamma (t) | \geq 2^{- N}$. \par

Let $1 \leq k \leq N$ such that $2^{- k} \leq 1 - |\gamma (t) | < 2^{- k + 1}$. Since $|p_{j, k}| = 1 - 2^{- k}$, we have:
\begin{displaymath} 
|\gamma (t) - p_{j, k} | \leq |\gamma (t) - \xi_j| + |\xi_j - p_{j, k}| \leq C (1 - |\gamma (t)|) + 2^{- k} \leq (2 C + 1) 2^{ - k} .
\end{displaymath} 
Hence 
\begin{displaymath} 
| \gamma (t) - p_{j,k} | \leq M \min \big(1 - | \gamma (t) | , 1 - | p_{j, k} | \big) \,,
\end{displaymath} 
with $ M = 2 C + 1$. By \eqref{petit lemme}, we get  
$\Big|\frac{\gamma (t) - p_{j k}} {1 - \overline{p_{j, k}}\, \gamma (t)} \Big| \leq \chi$, where $\chi = M / \sqrt{M^2 + 1}$ is $< 1$, and 
therefore $|B [\gamma (t)] | \leq \chi^d$. 
\qed
\par
\medskip

We can now end the proof of Theorem~\ref{general}. We get: 
\begin{align*} 
\frac{1}{h_n} \int_{ \{ s_N < |t - t_j| \leq s_n \} \cap \{ |t - t_j| \leq r_j\}} \hskip -80 pt 
&  \hskip 80 pt| B \big( \gamma (t) \big) |^2 \, \frac{dt}{2 \pi} \\
& \leq \frac{1}{h_n}  \, \frac{2 s_n}{2 \pi} \, \chi^{2 d} 
=  \frac{1}{h_n}  \, \frac{\omega^{- 1} (c^{- 1} C \, 2^{- n})}{\pi} \, \chi^{2 d} \\ 
& =  \frac{c^{- 1} C}{\pi} \, \frac{\omega^{- 1} (c^{- 1} C \, 2^{- n})}{c^{- 1} C \, 2^{- n}} \, \chi^{2 d} \\ 
& \leq \frac{1}{\pi} \, \omega^{- 1} (c^{- 1} C) \, \chi^{2 d} ,
\end{align*} 
since $\omega^{- 1} (x)/ x$ is non-decreasing. \par
\smallskip

We therefore get, setting $\kappa = c^{- 1} C$ and $L = L_1 + \cdots +L_p$:
\begin{align*}
\frac{1}{h_n} \sum_{j = 1}^p L_j \int_{\{| \gamma (t) - \xi_j | \leq C h_n \}\cap \{|t - t_j| \leq r_j\}} \hskip -80 pt 
& \hskip 80 pt |B [\gamma (t)]|^2 \, \frac{dt}{2 \pi} \\
& \leq \frac{\kappa L}{\pi} \, \frac{\omega^{- 1} (\kappa \, 2^{- N})}{\kappa \, 2^{- N}} 
+ \frac{L \, \omega^{- 1} (\kappa)}{\pi} \chi^{2 d} . 
\end{align*}

Choose now $d = d_N$, where $d_N$ is defined by \eqref{d_N}, with $\sigma = 1/ \log ( \chi^{- 2})$. Then 
$\chi^{2d} \leq \omega^{- 1} (\kappa \, 2^{- N}) /(\kappa \, 2^{- N})$, and, since the Blaschke product $B$ has now $p \, N d_N$ zeroes, we get, for some 
positive constant $K$:
\begin{displaymath} 
a_{p N d_N + 1} (C_\phi) \leq K \, \sqrt{\frac{\omega^{ - 1}  (\kappa \, 2^{- N})}{\kappa \, 2^{- N}}} \, \raise 1pt \hbox{,}
\end{displaymath} 
and that ends the proof of Theorem~\ref{general}. 
\qed 
\bigskip  

\noindent{\bf Proof of Theorem~\ref{polygon}.} It suffices to consider the case when $\phi$ is a conformal map from $\D$ onto ${\mathbf P}$. Indeed, let 
$\psi$ be such a conformal map. In the general case, our assumption allows to write $\phi = \psi \circ u$, where $u = \psi^{- 1}\circ \phi \colon \D\to \D$ is 
analytic. It follows that $C_\phi = C_u \circ C_\psi$ and that $a_n (C_\phi) \leq \| C_u \| \, a_n (C_\psi)$. Therefore, we may and shall assume that $\phi$ 
itself is this conformal map. \par
\smallskip

Let us denote by $\xi_1, \ldots, \xi_p$ the vertices of ${\mathbf P}$. Let $0 < \pi \mu_j < \pi$ be the exterior angle of ${\mathbf P}$ at $\xi_j$, namely the 
complement to $\pi$ of the interior angle; so that: 
\begin{displaymath} 
\sum_{j = 1}^p\mu_j = 2 \,, \quad \text{and} \quad  0 < \mu_j < 1.
\end{displaymath} 
If one sets $\theta_j = 1 - \mu_j$, one has $0 < \theta_j < 1$. \par
\smallskip

We then use the explicit form of $\phi$ given by the Schwarz-Christoffel formula (\cite{Nehari}, page 193):  
\begin{equation} \label{three} 
\phi (z) = A \int_0^z \frac{d w}{(a_1 - w)^{\mu_1} \cdots (a_p - w)^{\mu_p}} + B \,,
\end{equation}
for some constants $A \neq 0$ and $B \in \C$ and where $a_1, \ldots, a_p \in \partial \D$ are such that $\xi_j = \phi (a_j)$, $j = 1, \ldots, p$. 
If, as before, we write $\gamma (t) = \phi (\e^{it})$, we have $\xi_j = \gamma (t_j)$, with $a_j = \e^{i t_j}$ (note that here $E_{\xi_j} = \{ t_j\}$). \par

As we already said, condition \eqref{reluc glob} is trivially satisfied for a polygon. \par 

To end the proof, we use Theorem~\ref{general} and its Example 1. For that it suffices to  show that, for $| t - t_j |$ small enough, we have:
\begin{equation} \label{six}  
| \gamma (t) - \xi_j |  \approx | t - t_j |^{\theta_j} . 
\end{equation}
If $z\in \D$ is close to $a_j$, it follows from \eqref{three} that we can write
\begin{displaymath}
\phi (z) = A \int_0^z f_j (w) \frac{dw}{(a_j - w)^{\mu_j}} + B , 
\end{displaymath}
where $f_j$ is holomorphic near $a_j$ and $f_j (a_j) \neq 0$ since 
\begin{displaymath}
| f_j (a_j) | =\prod_{k\neq j, 1\leq k\leq p} | a_j - a_k |^{- \mu_k} . 
\end{displaymath}
Write $f_j (w) = f_j (a_j) + (a_j - w) g_j (w)$ where $g_j$ is holomorphic near $a_j$. We get:
\begin{align*}
\phi (z) 
& = A f_j (a_j) \int_0^z \frac{dw}{(a_j - w)^{\mu_j}} +  B + \int_0^z g_j (w) (a_j - w)^{\theta_j} \, dw \\
& := A f_j (a_j) \int_0^z \frac{dw}{(a_j - w)^{\mu_j}} + B +\psi_j (z), 
\end{align*}
which can still be written (since $\theta_j > 0$):
\begin{equation} \label{debloque} 
\phi (z) = \lambda_j (a_j - z)^{\theta_j} + c_j + \psi_j (z) ,
\end{equation}
where $\lambda_j \neq 0$, $c_j\in \C$, $\psi_j $ is Lipschitz near $a_j$ and $\xi_j = \phi (a_j) = c_j + \psi_j (a_j)$. 
Now, we easily get \eqref{six}. Indeed, for $t$ near $t_j$, it follows from \eqref{debloque} that (recall that $\gamma (t) = \phi (\e^{it})$ and 
$\gamma (t_j) = \xi_j$):
\begin{displaymath}
| \gamma (t) - \gamma (t_j) | = | \lambda_j | \, | \e^{it} - \e^{it_j} |^{\theta_j} + O \,  (| t - t_j |) ,
\end{displaymath}
which the claimed estimate \eqref{six} since $\lambda_j \neq 0$ and $| t - t_j |$ is negligible compared to  
$| t - t_j |^{\theta_j} \approx | \e^{it} - \e^{it_j} |^{\theta_j}$. 
\qed


\section{Lower bound and radial behavior} \label{lower} 

We shall consider symbols $\phi$ taking real values in the real axis ({\it i.e.} its Taylor series has real coefficients) and such that 
$\lim_{r \to 1^{-}} \phi (r) = 1$, with a given speed. 

\begin{definition} 
We say that the analytic map $\phi \colon \D \to \D$ is \emph{real} if it takes real values on $]- 1, 1[$, and that $\phi$ is an 
\emph{$\omega$-radial symbol} if it is real and there is a modulus of continuity $\omega \colon [0, 1] \to [0, 2]$ such that:
\begin{equation} \label{speedy} 
\qquad \qquad 1 - \phi (r) \leq \omega (1 - r) \,, \qquad 0 \leq r < 1 \,.
\end{equation} 
\end{definition}

With those  definitions and notations, one has:

\begin{theorem} \label{lastmin} 
Let $\phi$ be a real and $\omega$-radial symbol. Then, for the approximation numbers $a_n (C_\phi)$ of the composition operator $C_\phi$ of symbol $\phi$, 
one has the following lower bound:
\begin{equation} \label{general lower bound} 
a_n (C_\phi) \geq  c \, \sup_{0 < \sigma < 1} \sqrt{\frac{\omega^{- 1}(a \, \sigma^n)}{a\, \sigma^n}} \exp\bigg[- \frac{20}{1 - \sigma}\bigg] , 
\end{equation}
where $a = 1 - \phi (0) > 0$ and $c$ is another constant depending only on $\phi$. 
\end{theorem} 

Observe that, for the lens map $\lambda_\theta$ (see \cite{lens}, Lemma~2.5), we have $\omega^{- 1} (h) \approx h^{1/\theta}$, so that adjusting 
$\sigma = 1 - 1/\sqrt n$, we get 
\begin{equation} \label{mino lens}
a_n (C_{\lambda_\theta}) \geq c \, \exp \big( - C \sqrt{n} \big) , 
\end{equation} 
which is the result of \cite{LIQUEROD}, Proposition~6.3. \par

For the cusp map $\phi$ (see Section~\ref{section cusp}), we have $\omega^{- 1} (h) \approx \e^{- C' / h}$, so that taking $\sigma = \exp (- \log n / 2n)$, 
we get:
\begin{equation} \label{mino cusp}
a_n (C_\phi) \geq c \exp (- C \, n / \log n ) .
\end{equation} 

We shall use the same methods as for lens maps (see \cite{LIQUEROD}, Proposition~6.3). \par \smallskip

We need a lemma. Recall (see \cite{Hoffman-livre} pages~194--195, or \cite{Nikolski-livre} pages~302--303) that if $(z_j)$ is a Blaschke sequence, its Carleson 
constant $\delta$ is defined as $\delta = \inf_{j \geq 1} (1 - |z_j|^2) \, |B' (z_j)|$, where $B$ is the Blaschke product whose zeros are the $z_j$'s. Now 
(see \cite{Garnett}, Chapter~VII, Theorem~1.1), every $H^\infty$-interpolation sequence $(z_j)$ is a Blaschke sequence and its Carleson constant 
$\delta$ is connected to its interpolation constant $C$ by the inequalities
\begin{equation} 
1/ \delta \leq C \leq \kappa / \delta^2
\end{equation} 
where $\kappa$ is an absolute constant (actually $C \leq \kappa_1 (1 / \delta) ( 1 + \log 1/ \delta) $). Now, if $(z_j)$ is a $H^\infty$-interpolation 
sequence with constant $C$, the sequence of the normalized reproducing kernels  $f_j =  K_{z_j} / \|K_{z_j} \|$ satisfies 
\begin{displaymath} 
C^{- 1} \big( \sum | \lambda_j|^2 \big)^{1/2} \leq \big\| \sum \lambda_j f_j \|_{H^2} \leq C \, \big( \sum | \lambda_j|^2 \big)^{1/2} 
\end{displaymath} 
(see \cite{LIQUEROD}, Lemma~2.2). 

\begin{lemma} \label{lemme avec Carleson} 
Let $\phi \colon \D \to \D$ be an analytic self-map. Let $u = (u_1, \ldots, u_n)$ be a finite sequence in $\D$ and set $v_j = \phi (u_j)$, $v = (v_1, \ldots, v_n)$. 
Denote by $\delta_v$ the Carleson constant of the finite sequence $v$ and set 
\begin{displaymath} 
\mu_n^2 = \inf_{1 \leq j \leq n} \frac{1 - | u_j |^2}{1 - | \phi (u_j) |^2} \, \cdot
\end{displaymath} 
Then, for some constant $c' > 0$, we have the lower bound:
\begin{equation} \label{aprio} 
a_n (C_\phi) \geq c' \, \delta_{v}^4 \, \mu_n . 
\end{equation}
\end{lemma} 

\noindent{\bf Proof.} Recall first that the Carleson constant $\delta$ of a Blaschke sequence $(z_j)$ is also equal to:
\begin{displaymath} 
\delta = \inf_{k\geq 1} \prod_{j\neq k} \rho (z_k, z_j) \,,
\end{displaymath} 
where $\rho (z, \zeta) = \big| \frac{z - \zeta}{1 - \overline{z} \, \zeta} \big|$ is the pseudo-hyperbolic distance between $z$ and $\zeta$. Now, the 
Schwarz-Pick Lemma (see \cite{Beardon-Minda}, Theorem~3.2) asserts that every analytic self-map of $\D$ contracts the pseudo-hyperbolic distance. 
Hence $\rho \big( \phi (u_j), \phi (u_k) \big) \leq \rho (u_j, u_k)$ and so, if $\delta_u$ and $\delta_v$ denote the Carleson constants of $u$ and $v$:
\begin{displaymath} 
\delta_u \geq \delta_v. 
\end{displaymath} 
\par\smallskip

Let now $R$ be an operator of rank $< n$. There exists a function $f = \sum_{j=1}^n \lambda_j K_{u_j} \in H^2 \cap \ker R$ with 
$\| f \| =1$. We thus have: 
\begin{align*}
\| C_\phi ^* - R \|^2 
& \geq \| C_\phi ^* (f) - R (f) \|_2^2 = \| C_\phi ^* (f) \|_2^2  
= \bigg\| \sum_{j=1}^n \lambda_j K_{v_j} \bigg\|_2^2 \\
& \geq C_{v}^{- 2} \sum_{j = 1}^n | \lambda_j |^2 \| K_{v_j} \|_2^2
= C_{v}^{- 2} \sum_{j = 1}^n \frac{| \lambda_{j} |^2}{1 - |v_j|^2} \\
& \geq C_{v}^{- 2} \,  \mu_n^2 \sum_{j=1}^n \frac{| \lambda_j |^2}{1 - | u_j |^2} \\
& \geq C_{u}^{- 2} C_{v}^{- 2} \, \mu_n^2  \| f \|_2^2 
= C_{u}^{- 2} C_{v}^{- 2} \, \mu_n^2 \\
& \geq \kappa^{- 4} \, \delta_u^4 \, \delta_v^4 \, \mu_n^2 
\geq \kappa^{- 4} \, \delta_v^8 \, \mu_n^2 ,
\end{align*}
and hence $a_n (C_\phi) \geq \kappa^{- 2} \, \delta_v^4 \, \mu_n$. \qed
\bigskip 

\noindent{\bf Remark.} This lemma allows to give, in the Hardy case, a simpler proof of Theorem~4.1 in \cite{LIQUEROD}, avoiding the use of Lemma~2.3 
and Lemma~2.4 (concerning the backward shift) in that paper. Recall that this theorem says that for every non-increasing sequence $(\eps_n)_{n \geq 1}$ of  
positive real numbers tending to $0$, there exists a univalent symbol $\phi$ such that $\phi (0) = 0$ and $C_\phi \colon H^2 \to H^2$ is compact, but 
$a_n (C_\phi) \gtrsim \eps_n$ for every $n \geq 1$. Let us sketch briefly the argument. We use the notation of \cite{LIQUEROD}, Lemma~4.6. The symbol 
$\phi$ is defined as $\phi (z) = \sigma^{- 1} (\e^{- 1} \sigma (z))$, where $\sigma$ is some conformal map $\sigma \colon \D \to \Omega$. We set 
$A_j = (1 /C_0) \log (1/ \eps_{j + 1})$, $r_j = \sigma^{- 1} (\e^j)$. Then $\phi (r_{j + 1}) = r_j$ and (see \cite{LIQUEROD}, pages 444--446):
\begin{displaymath} 
\frac{1 - r_{j + 1}}{1 - r_j} \geq \exp (- 2 C_0 A_j) .
\end{displaymath} 
We shall apply the above Lemma~\ref{lemme avec Carleson} with $u_j = r_j $. Then $v_j = \phi (u_j) = r_{j - 1}$. Hence
\begin{displaymath} 
\frac{1 - |u_j|^2}{1 - |v_j|^2} \geq \frac{1}{2}\, \frac{1 - |u_j|}{1 - |v_j|}  = \frac{1}{2}\, \frac{1 - r_j }{1 - r_{j - 1}}  
\geq \frac{1}{2}\, \exp (-2 C_0 A_{j - 1}) = \frac{1}{2} \, \eps_j^2 \geq \frac{1}{2} \, \eps_n^2 .
\end{displaymath} 
It follows that $\mu_n \geq \eps_n / \sqrt{2}$. \\
On the other hand, $(r_j)_{j \geq 1}$ is an interpolating sequence (see \cite{LIQUEROD}, Lemma~4.6); hence there is a constant $\delta > 0$ (which does not 
depend on $n \geq 1$) such that $\delta_v \geq \delta$. Therefore Lemma~\ref{lemme avec Carleson} gives 
\begin{displaymath} 
a_n (C_\phi) \geq c\, \delta^4 \eps_n \,,
\end{displaymath} 
which gives Theorem~4.1 of \cite{LIQUEROD}. \qed

\bigskip

\noindent{\bf Proof of Theorem~\ref{lastmin}.} Fix $0 < \sigma < 1$ and define inductively $u_j \in [0, 1)$ by $u_0 = 0$ and the relation
\begin{displaymath} 
\qquad 1 - \phi (u_{j+1}) = \sigma [1 - \phi (u_j) ] \quad \text{with } 1 > u_{j + 1} >  u_j  
\end{displaymath} 
(using the intermediate value theorem). \par

Setting $v_j = \phi (u_j)$, we have $-1 < v_j < 1$, 
\begin{equation} \label{catal} 
\frac{\ 1 -  v_{j + 1} }{1 - v_j \ } = \sigma, 
\end{equation}
and 
\begin{equation} \label{catal bis}
\qquad \qquad \qquad  1 - v_n = a \, \sigma^n \, , \qquad \quad \text{with } a = 1 - \phi (0) . \qquad  
\end{equation} 
Now observe that, for $1 \leq j \leq n$, one has, due to the positivity of $u_j$ and $v_j$, to \eqref{speedy}, and the fact that 
$r_\omega (x) = \omega^{- 1} (x) / x$ is increasing: 
\begin{align*} 
\frac{1 - | u_j |^2}{1 - | v_j |^2} 
& \geq \frac{1 - u_j}{2 (1 - v_j)} \geq \frac{1}{2} \, \frac{\omega^{- 1} (1 - v_j)}{1 - v_j} = \frac{1}{2}\, r_\omega (1 - v_j) \\
&  \geq \frac{1}{2}\, r_\omega (1 - v_n) = \frac{1}{2}\, r_\omega (a \, \sigma^n), 
\end{align*} 
which proves that $\mu_n^2 \geq r_\omega (a \, \sigma^n) / 2$. Furthermore, the sequence $(v_j)$ satisfies, by \eqref{catal}, a condition very similar to 
Newman's condition with parameter $\sigma$. In fact, for $k > j$, we have 
\begin{displaymath} 
\frac{|v_k - v_j|}{|1 - v_k v_j|} = \frac{(1 -v_j) - (1 - v_k)}{(1 - v_j) + v_j (1 - v_k)} \geq \frac{(1 -v_j) - (1 - v_k)}{(1 - v_j) + (1 - v_k)} 
= \frac{1 - \sigma^{k - j}}{1 + \sigma^{k - j}} \, \cdot
\end{displaymath} 
Analogously, for $j > k$, we have $\frac{|v_k - v_j|}{|1 - v_k v_j|} \geq \frac{1 - \sigma^{j - k}}{1 + \sigma^{j - k}}\,$. Thus, as in the proof of 
\cite{DUREN-livre}, Theorem~9.2, we have, for every $k$, 
\begin{displaymath} 
\prod_{j \neq k} \rho (v_j, v_k) = \prod_{j \neq k} \frac{|v_k - v_j|}{|1 - v_k v_j|} 
\geq \prod_{l = 1}^\infty \Big( \frac{1 - \sigma^l}{1 + \sigma^l} \Big)^2 \, .
\end{displaymath} 
Consequently, $\delta_v \geq \prod_{l = 1}^\infty \big( \frac{1 - \sigma^l}{1 + \sigma^l} \big)^2 \geq \exp \big(- \frac{5}{1 - \sigma} \big)$, by 
\cite{LIQUEROD}, Lemma~6.4. Finally, use \eqref{aprio} to get: 
\begin{displaymath} 
a_n (C_\phi) \geq c' \, \delta_v^4 \, \mu_n \geq c \, \exp \Big(- \frac{20}{1 - \sigma} \Big) \sqrt{r_\omega (a \, \sigma^n)}. 
\end{displaymath} 
Taking the supremum over $\sigma$, that ends the proof of Theorem~\ref{lastmin}. \qed
\bigskip \goodbreak 

\noindent{\bf Remark.} The proof shows that 
\begin{equation} 
a_n (C_\phi) \geq \sup_{u_1, \ldots, u_n \in (0, 1)} \inf _{
\substack{\scriptscriptstyle f \in \langle K_{u_1}, \ldots, K_{u_n} \rangle \\ \scriptscriptstyle \| f  \| = 1}} \| C_\phi^\ast f \| \,,
\end{equation} 
where $\langle K_{u_1}, \ldots, K_{u_n} \rangle$ is the linear space generated by $n$ distinct reproducing kernels $K_{u_1}, \ldots, K_{u_n}$. But if $B$ is the 
Blaschke product with zeros $u_1, \ldots, u_n$, then $\langle K_{u_1}, \ldots, K_{u_n} \rangle = (B H^2)^\perp$, the \emph{model space} associated to $B$. 
Hence 
\begin{equation} 
a_n (C_\phi) \geq \sup_B  \inf _{
\substack{\scriptscriptstyle f \in (B \! H^2)^\perp \\  \scriptscriptstyle \| f  \| = 1}} \| C_\phi^\ast f \| \,,
\end{equation} 
where the supremum is taken over all Blaschke products with $n$ zeros on the real axis $(0, 1)$. This has to be compared with the upper bound (which gives 
\eqref{majo avec Blaschke}, see \cite{lens}, proof of Lemma~2.4):
\begin{equation} 
a_n (C_\phi) \leq \inf_B \big\| {C_\phi}_{\mid B H^2} \big\| 
= \inf_B \sup_{\substack{\scriptscriptstyle f \in B \! H^2 \\ \scriptscriptstyle \| f \| = 1}} \| C_\phi f \| \,,
\end{equation} 
where the infimum is over the Blaschke products with less than $n$ zeros (in the Hilbert space $H^2$, the approximation number $a_n (C_\phi)$ is equal to the 
Gelfand number $c_n (C_\phi)$, which is, by definition, less or equal to $\big\| {C_\phi}_{\mid B H^2} \big\|$, since $BH^2$ is of codimension $< n$). 

\section{Examples} 

\subsection{The cusp map} \label{section cusp}  

\begin{figure}[ht]
\centering
\includegraphics[width=6cm]{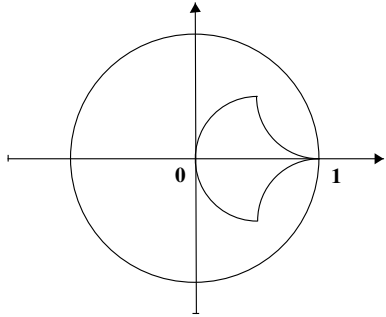}
\caption{\it Cusp map domain} \label{cusp map}
\end{figure}

\begin{definition}
The \emph{cusp map} is the conformal mapping $\phi$ sending the unit disk $\D$ onto the domain represented on Figure~\ref{cusp map}. 
\end{definition}

This map was first introduced in \cite{TRANS} (see also \cite{RACSAM}). Explicitly, $\phi$ is defined as follows. \par
We first map $\D$ onto the half-disk $\D^+ =\{z \in \D \, ; \ \Re z > 0\}$. To do that, map $\D$ onto itself by $z \mapsto i z$; then map  $\D$ onto the 
upper half-plane $\H = \{z \in \C \, ; \ \Im z > 0 \}$ by:
\begin{displaymath} 
T (u) = i \, \frac{1 + u}{1 - u} \,\cdot
\end{displaymath} 
Take the square root to map $\H$ in the first quadrant $Q_1 = \{z \in \H \, ; \ \Re z > 0\}$, and go back to the half-disk 
$\{z \in \D \, ; \ \Im z < 0\}$ by $T^{- 1}$: $T^{- 1} (s) = \frac{1 + is}{is - 1}$; finally, make a rotation by $i$ to go onto $\D^+$. We get:
\begin{equation} \label{phi zero}
\phi_0 (z) = \frac{\displaystyle \Big( \frac{z - i}{i z - 1} \Big)^{1/2} - i} {\displaystyle - i \, \Big( \frac{z - i}{i z - 1} \Big)^{1/2} + 1} \, \cdot
\end{equation} 
One has $\phi_0 (1) = 0$, $\phi_0 (- 1) = 1$, $\phi_0 (i) = - i$ and $\phi_0 (- i) = i$. The half-circle $\{z \in \T \, ; \ \Re z \geq 0 \}$ is mapped onto the 
segment $[- i , i ]$ and the segment $[- 1, 1]$ onto the segment $[0, 1]$. \par
Set now, successively, 
\begin{equation} 
\phi_1 (z) = \log \phi_0 (z), \quad \phi_2 (z) = - \frac{2}{\pi}\, \phi_1 (z) + 1, \quad \phi_3 (z) = \frac{1}{\phi_2 (z)}  \, \raise 1pt \hbox{,} 
\end{equation} 
and finally:
\begin{equation} 
\phi (z) = 1 - \phi_3 (z) \,.
\end{equation} 
Hence:
\begin{equation}
1 - \phi (z) = \frac{1} {1 + \frac{2}{\pi} \, \log \big( 1 / |\phi_0 (z) | \big) - i \frac{2}{\pi} \arg \phi_0 (z)} \, \cdot
\end{equation}
$\phi_2$ maps $\D$ onto the semiband $\{z\in \C \, ; \ \Re z > 1 \text{ and } |\Im z  | < 1\}$. One has $\phi (1) = 1$, $\phi (- 1) = 0$, 
$\phi (i) = (1 + i)/ 2$ and $\phi (- i) = (1 - i)/2$. \par 
The domain $\phi (\D)$ is edged by three circular arcs of radii $1/2$ and of respective centers $1/2$, $1 + i/2$ and $1 -i/2$. The real interval $] - 1, 1[$ 
is mapped onto the real interval $]\, 0, 1 [$ and the half-circle $\{\e^{i \theta} \, ; \ |\theta | \leq \pi/ 2\}$ is sent onto the two circular arcs tangent at $1$ 
to the real axis. 

\begin{lemma} \label{comportement local} \hfill \par
1) For $0 < r < 1$, let $\gamma = \frac{\pi}{4} - \arctan r = \arctan [ (1 - r) / (1 + r) ]$; then:
\begin{equation} \label{phi zero sur axe reel}
\phi_0 (r) = \tan (\gamma / 2) \,. 
\end{equation} 
Hence, when $r$ tends to $1_-$, one has:
\begin{equation} \label{phi sur axe reel}
1 - \phi (r) \sim \frac{\pi}{2} \, \frac{1}{\log (1 / \gamma)} \sim \frac{\pi}{2} \, \frac{1}{\log (1 / (1 - r))} \, \cdot
\end{equation} 
2) For $|\theta | < \pi / 2$, one has:
\begin{equation} \label{phi zero sur cercle}
\phi_0 (\e^{i \theta}) = - i\, \frac{\tan (\theta / 2)}{1 + \sqrt{1 - \tan^2 (\theta / 2)}} \,\cdot
\end{equation} 
Hence, when $\theta$ tends to $0$, one has:
\begin{equation} \label{phi sur cercle}
1 - \phi (\e^{i \theta}) \sim \frac{\pi}{2} \, \frac{1}{\log (1/ |\theta|)} \, \cdot
\end{equation} 
\end{lemma}
\noindent{\bf Proof.} 1) One has:
\begin{displaymath} 
T (i r) = \frac{r - i}{i r - 1} = - \frac{2r}{1 + r^2} + i\, \frac{1 - r^2}{1 + r^2} = - \sin \alpha + i \cos \alpha \,,
\end{displaymath} 
with $r = \tan (\alpha/2)$; hence $T (i r) = \cos (\alpha + \pi/2) + i \sin (\alpha + \pi / 2) = \e^{i (\alpha + \pi / 2)}$. Set 
$\beta = \frac{\alpha}{2} + \frac{\pi}{4}$; one gets:
\begin{displaymath} 
\phi_0 (r) = \frac{\e^{i \beta} - i}{- i \e^{i \beta} + 1} = \frac{\cos \beta}{1 + \sin \beta} 
= \frac{\sin \gamma}{1 + \cos \gamma} = \tan (\gamma / 2) \, 
\end{displaymath} 
with $\gamma = (\pi / 2) - \beta = (\pi/ 4) - (\alpha/ 2) = (\pi / 4) - \tan^{-1} r$. Then \eqref{phi sur axe reel} follows. \par
2) Let $\tau = \frac{\pi}{2} - \theta$; one has:
\begin{displaymath} 
T (i \e^{i \theta}) = \frac{\e^{i \theta} - i}{i \e^{i \theta} - 1} = \frac{- \cos \theta}{ 1 + \sin \theta} = \frac{- \sin \tau}{1 + \cos \tau} 
= - \tan (\tau / 2). 
\end{displaymath} 
Note that $0 < \tau / 2 < \pi / 2$ since $|\theta | < \pi /2$; hence $\tan (\tau / 2) > 0$. Therefore: 
\begin{displaymath} 
\phi_0 (\e^{i \theta}) = \frac{i \sqrt{\tan (\tau/ 2)} - i}{- i . i \sqrt{\tan (\tau / 2)} + 1} 
= i\, \frac{ \sqrt{\tan (\tau/ 2)} - 1}{ \sqrt{\tan (\tau/ 2)} + 1} \, \cdot
\end{displaymath} 
But 
\begin{displaymath} 
\tan (\tau / 2) = \tan \Big(\frac{\pi}{4} - \frac{\theta}{2} \Big) = 
\frac{1 - \tan (\theta / 2)}{1 + \tan (\theta / 2)} \,;
\end{displaymath} 
it follows that:
\begin{align*} 
\phi_0 (\e^{i \theta}) 
& = i\,  \frac{\sqrt{1 - \tan (\theta / 2) } - \sqrt{1 + \tan (\theta / 2)}} {\sqrt{1 - \tan (\theta / 2) } + \sqrt{1 + \tan (\theta / 2)}} \\
& = i \, \frac{\big( 1 - \tan (\theta / 2) \big) - \big( 1 + \tan (\theta / 2) \big)} {\big( \sqrt{1 - \tan (\theta / 2) } + \sqrt{1 + \tan (\theta / 2)} \big)^2} \\ 
& = - i\, \frac{\tan (\theta / 2)} {1 + \sqrt{1 - \tan^2 (\theta / 2)} } \, \cdot
\end{align*} 

Now, since $\phi_0 (\e^{i \theta}) \sim - i \theta/ 4$ as $\theta$ tends to $0$, we get \eqref{phi sur cercle}. \qed
\par\bigskip

It follows from this lemma and from Theorem~\ref{general} and Theorem~\ref{lastmin} that one has the following estimate. 
\begin{theorem} \label{estimation number}
For the approximation numbers $a_n (C_\phi)$ of the composition operator $C_\phi \colon H^2 \to H^2$ of symbol the cusp map $\phi$, we have:
\begin{equation} 
\qquad \qquad \qquad \qquad \e^{- c_1 \, n / \log n} \lesssim a_n (C_\phi) \lesssim \e^{- c_2 \, n / \log n} , \qquad \quad n = 2, 3, \ldots, 
\end{equation} 
for some constants $c_1 > c_2 > 0$.
\end{theorem} 

\noindent{\bf Proof.} 1) \emph{Upper estimate}. Note first that, since the domain $\phi (\D)$ is contained in the right half-plane and in the symmetric angular 
sector of vertex $1$ and opening $\pi/ 2$, there is a constant $C > 0$ such that $|1 - \gamma (t) | \leq C\, (1 -|\gamma (t)|)$ and we have 
\eqref{reluc}. Then \eqref{phi sur cercle} in Lemma~\ref{comportement local} gives \eqref{ponct}. The upper estimate is hence given in 
Theorem~\ref{general}  and \eqref{majo cusp}. \par
2) \emph{Lower estimate}. By Lemma~\ref{comportement local}, \eqref{phi sur axe reel}, one has \eqref{speedy}. Since $\phi$ is a real symbol, the upper 
estimate follows from Theorem~\ref{lastmin}, and \eqref{mino cusp}. \qed

\subsection{The Shapiro-Taylor map} \label{ST map}

This one-parameter map $\varsigma_\theta\,$, $\theta > 0$, was introduced by J. Shapiro and P. Taylor in 1973 (\cite{SHTA}) and was further studied, with a 
slightly different definition, in \cite{JFA}, Section~5. J. Shapiro and P. Taylor proved that $C_{\varsigma_\theta} \colon H^2 \to H^2$ is always compact, but is 
Hilbert-Schmidt if and only if $\theta > 2$. It is proved in \cite{JFA}, Theorem~5.1, that $C_{\varsigma_\theta}$ is in the Schatten class $S_p$ if and only if 
$p > 4 / \theta $. 
\par\smallskip 

Here, we shall use these maps $\varsigma_\theta$ to see the limitations of our previous methods. \par
\medskip

We first recall their definition. \par\smallskip

For $\eps > 0$, we set $V_\eps = \{ z\in \C \, ; \ \Re z > 0 \text{ and } |z | < \eps \}$. For $\eps = \eps_\theta > 0$ small enough, one can define 
\begin{equation} 
f_\theta (z) = z (- \log z )^\theta ,
\end{equation} 
for $z \in V_\eps$, where $\log z$ will be the principal determination of the logarithm. Let now $g_\theta$ be the conformal mapping from $\D$ onto $V_\eps$, 
which maps $\T = \partial \D$ onto $\partial V_\eps$, defined by $g_\theta (z) = \eps\, \phi_0 (z)$, where $\phi_0$ is given in \eqref{phi zero}. \par
Then, we define:
\begin{equation} 
\varsigma_\theta = \exp ( - f_\theta \circ g_\theta) .
\end{equation} 
One has $\varsigma_\theta (1) = 1$ and $g_\theta (\e^{it}) \sim - i t /4$ as $t$ tends to $0$, by Lemma~\ref{comportement local}; hence, when $t$ is near 
of $0$: 
\begin{displaymath} 
| 1 - \varsigma_\theta (\e^{it}) | \approx |f_\theta [ g_\theta (\e^{it}) ]| \approx |t| \, [\log (1 / |t|) ]^\theta .
\end{displaymath} 

If we were allowed to apply Theorem~\ref{general}, we would get that $a_n (C_{\varsigma_\theta}) \lesssim 1/ n^{\theta / 4}$, which would be in accordance 
with the fact that $C_{\varsigma_\theta}$ is in the Schatten class $S_p$ if and only if $p > 4 / \theta$. However, condition \eqref{reluc} is not satisfied: by 
\cite{JFA}, equations (5.5) and (5.6), one has $1 - |\varsigma_\theta (\e^{it}) | \approx |t| (\log 1/ |t|) ^{\theta - 1}$, whereas 
$| 1 - \varsigma_\theta (\e^{it}) | \approx |t| (\log 1 / |t|)^\theta$. \par
\medskip

On the other hand, by the Lemma~\ref{comportement local} again, $g_\theta (r)  \sim \eps (1 - r)/4$ as $r$ tends to $1$; hence, when $r$ is near to $1$:
\begin{displaymath} 
1 - \varsigma_\theta (r) \approx (1 - r) \big( \log 1 / (1 - r) \big)^\theta ,
\end{displaymath} 
so $\varsigma_\theta$ is a real $\omega$-radial symbol with $\omega (t) = t (\log 1 / t)^\theta$. Hence, we get from Theorem~\ref{lastmin}:  
\begin{displaymath} 
a_n (C_{\varsigma_\theta}) \gtrsim \frac{1}{n^{\theta /2}} \, \raise 1pt \hbox{,}
\end{displaymath} 
taking $\sigma = 1 / \e $ in \eqref{general lower bound}. However, this lower estimate is not the right one, since 
$C_{\varsigma_\theta}$ is in $S_p$ if and only if $p > 4 / \theta$. 


\section{Contact points} \label{application} 

It is well-known (and easy to prove) that for every compact composition operator $C_\phi \colon H^2 \to H^2$, the set of contact points 
\begin{displaymath} 
E_\phi = \{ \e^{i \theta} \, ; \ |\phi^\ast (\e^{i \theta} ) | = 1 \} 
\end{displaymath} 
has Lebesgue measure $0$. A natural question is: to what extent is this negligible set arbitrary? The following partial answer was given by E.A.~Gallardo-Guti\'errez 
and M.J.~Gonz\'alez in \cite{GALGON}. 
\begin{theorem} [E.A.~Gallardo-Guti\'errez and M.J.~Gonz\'alez] \label{calgon} 
There is a compact composition operator $C_\varphi$ on $H^2$  such that the Hausdorff dimension of $E_\phi$ is one.
\end{theorem} 

This was generalized by O. El-Fallah, K. Kellay, M. Shabankhah, and  H. Youssfi (\cite{EKSY}, Theorem~3.1):
\begin{theorem} [O. El-Fallah, K. Kellay, M. Shabankhah, H. Youssfi] \hfill \\ 
For every compact set $K$ of measure $0$ in $\T$, there exists a Schur function $\phi \in A (\D)$, the disk algebra, such that the associated composition 
operator $C_\phi$ is Hilbert-Schmidt on $H^2$ and $E_\phi = K$.
\end{theorem} 

As an application of our previous results, we shall extend these results, with a very simple proof. Our composition operator will not even be compact, or 
Hilbert-Schmidt, but in all Schatten classes $S_p$, and moreover its approximation numbers will be as small as possible. 

\begin{theorem} 
Let $K$ be a Lebesgue-negligible compact set of the circle $\T$. Then, there exists a Schur function $\psi \in A (\D)$, the disk algebra, such that $E_\psi = K$, 
$\psi (\e^{i \theta}) = 1$ for all $\e^{i \theta} \in K$,  and:  
\begin{equation} 
a_n (C_\psi) \leq a \exp(- b \, n / \log n) .
\end{equation} 
In particular, $C_\psi \in \bigcap_{p > 0} S_p$.
\end{theorem}

\noindent{\bf Proof.} According to the Rudin-Carleson theorem (\cite{BIS}), we can find $\chi \in A(\D)$ such that 
\begin{displaymath} 
\qquad \chi = 1 \text{ on } K \quad \text{and}  \quad |\chi | < 1 \text{ on } \overline{\D} \setminus K. 
\end{displaymath} 
Consider now the cusp map $\phi$, defined in Section~\ref{section cusp}. One has $\phi \in A (\D)$, $\phi (1) = 1$ and 
\begin{displaymath} 
a_{n} (C_\phi) \leq a' \exp (- b \, n / \log n).
\end{displaymath} 
We now spread the point $1$ by composing with the function $\chi$, which is equal to $1$ on the whole of $K$. 
We check that the composed map $\psi = \phi \circ \chi$ has the required properties. \par
That $\psi \in A (\D)$ is clear. For $z \in K$, one has $\psi (z) = \phi (1) = 1$, and for $z \in \overline{\D} \setminus K$, one has 
$| \chi (z) | < 1$; hence $| \psi (z) | < 1$. \par
To finish, since $C_\psi = C_\chi \circ C_\phi$, we have 
\begin{displaymath} 
a_{n} (C_\psi) \leq \| C_\chi \| \, a_{n} (C_\phi) \leq a'  \, \| C_\chi \| \, \exp(- b \, n / \log n) := \sigma_n,
\end{displaymath} 
proving the result (with $a = a' \, \| C_\chi \|$), since clearly  $\sum_{n = 1}^\infty \sigma_{n}^p < \infty$ for each $p > 0$.  
\qed
\bigskip
 
Actually, we can improve on the previous theorem by proving the following result. This result is optimal because if $\| \psi \|_\infty = 1$, we know  
(see \cite{LIQUEROD}, Theorem~3.4) that $\liminf_{n \to \infty} [a_n (C_\psi)]^{1/n} = 1$, so we cannot hope to get rid with the forthcoming vanishing 
sequence $(\eps_n)_n$.
  
\begin{theorem}\label{cherry} 
Let $K$ be a Lebesgue-negligible compact set of the circle $\T$ and $(\eps_n)_n$ a sequence of positive real numbers with limit zero. Then, there exists a Schur 
function $\phi\in A(\D)$ such that $E_\phi = K$, $\phi (\e^{i \theta}) = 1$ for all $\e^{i \theta} \in K$,  and 
\begin{equation} 
a_{n}(C_\phi)\leq C \exp (- n \,\eps_n) \,,
\end{equation} 
where $C$ is a positive constant.
\end{theorem}

This theorem is a straightforward consequence of  the following lemma. Recall that the Carleson function of the Schur function $\psi \colon \D \to \D$ is defined 
by: 
\begin{displaymath} 
\rho_\psi (h) = \sup_{|\xi | = 1} m (\{ t \in \T \, ; \ |\psi (\e^{i t})| \geq 1 - h \text{ and } | \arg ( \psi (\e^{i t}) \, \bar{\xi} ) | \leq \pi h\}). 
\end{displaymath} 

\begin{lemma} \label{end of the world} 
Let $\delta$ be a nondecreasing positive function on $(0, 1]$ tending to $0$  as $h \to 0$. Then, there exists a Schur function $\psi \in A(\D)$ such that 
$\psi (1) = 1$, $| \psi (\xi) |  < 1$ for $\xi \in \T \setminus \{1\}$, and such that $\rho_{\psi} (h) \leq \delta (h)$, for $h > 0$ small enough.
\end{lemma}

Once we have the lemma, in view of the upper bound in \cite{LIQUEROD}, Theorem~5.1, for approximation numbers, we can adjust the function 
$\delta$ so as to have $a_n (C_\psi) \leq K \e^{- n \eps_n}$. Then, we  compose $\psi$ with a peaking function $\chi$ as in the previous section and the map 
$\phi = \psi \circ \chi$ fulfills the requirements of Theorem~\ref{cherry}, with $C = K  \| \chi \|$.  \qed
\medskip

\begin{figure}[ht]
\centering
\includegraphics[width=6cm]{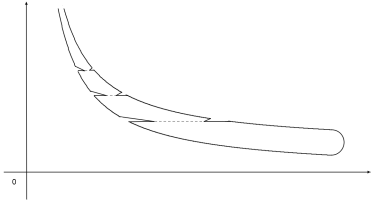}
\vskip - 10 pt 
\caption{\it Domain $\Omega$} \label{domain}
\end{figure}

\noindent{\bf Proof of Lemma~\ref{end of the world}.} We use a slight modification of the map $g$ constructed in \cite{REVI}, pages 66--67. Instead of 
taking a conformal map from $\D$ to the domain used in \cite{REVI}, we modify this domain by limiting it to the right- hand  side (by, say, a semicircle), as on the 
Figure~\ref{domain}. Let $\Omega$ this domain. This domain is limited by the two hyperbolas $y = 1/x$ and $y = (1/x) + 4\pi$. The limiting 
semicircle is chosen in order that $\Im w \geq 1$ for $w \in \Omega$. The lower part of the ``saw-teeth'' have an imaginary part equal to $4\pi n$. 
If $a \in \Omega$ is fixed and $\Omega_n$ is the part of the domain $\Omega$ such that $\Im w < 4 \pi n$, the horizontal sizes of the ``saw-teeth'' are 
chosen in order that the harmonic measure $\omega_\Omega (a, \partial \Omega \setminus \partial \Omega_n)$ is $\leq \delta_n := \delta (1/ 16 \pi (n + 1))$. 
Note that $\partial \Omega \setminus \partial \Omega_n \supseteq \{ w \in \partial \Omega \, ; \ \Im w > 4 \pi n \}$ (see \cite{REVI}, Lemma~4.2). 
\par

By Carath\'eodory-Osgood's Theorem (see \cite{Palka}, Theorem~IX.4.9), there is a unique homeomorphism $g$ from $\overline{\D}$ onto 
$\overline{\Omega} \cup\{\infty\}$ which maps conformally $\D$ onto $\Omega$ and such that $g (0) = a$ and $g (1) = \infty$ (we may choose these two 
values because if $h \colon \overline{\D} \to \overline{\Omega} \cup\{\infty\}$ is such a map, and $u$ is the automorphism of $\overline{\D}$ such that 
$u (0) = h^{- 1} (a)$ and $u (1) = h^{- 1} (\infty)$, then $g = h \circ u$ suits -- alternatively, having choosen $h (0) = a$, then, if $h (\e^{i \theta_0}) = \infty$, 
we take $g (z) = h (\e^{i \theta_0} z)$). \par

We define $\psi = (g - i) / (g + i)$. Then $\psi \colon \D \to \D$ is a Schur function and $\psi \in A (\D)$. Moreover, since the domain $\Omega$ is bounded 
horizontally, we have $\psi (1) = 1$ and $|\psi ( \e^{it}) | < 1$ for $0 < t < 2\pi$. \par

Now, $\rho_\psi (h) \leq m \big( \{ z \in \T \, ; \ | \psi (z)| > 1 - h \}\big) $. Writing $g = u  + i v$, one has:
\begin{displaymath} 
| \psi |^2  = \frac{u^2 + (v - 1)^2}{u^2 + (v + 1)^2} = 1 - \frac{4 v}{u^2 + (v + 1)^2} \, \cdot
\end{displaymath} 
Since $(1 - h)^2 \geq 1 - 2h$, the condition $| \psi (z)| > 1 - h$ implies that $\frac{2v}{u^2 + (v + 1)^2} \leq h$. But $0 < u \leq 1 + 2\pi \leq 8$ and 
$(v + 1)^2 \leq 4 v^2$ (since $v \geq 1$); we get hence $\frac{v}{32 + 2 v^2} \leq h$, or $\frac{32}{v} + 2v \geq \frac{1}{h}$. Using again the fact that 
$v \geq 1$, one obtains $2v \geq \frac{1}{h} - 32$, and hence $2v \geq \frac{1}{2h}$ for $0 < h \leq 1/64$. Therefore, for $0 < h \leq 1/64$, 
\begin{displaymath} 
\rho_\psi (h) \leq m (\{ z \in \T \, ; \ \Im \psi (z) \geq 1 /4h \}) \,.
\end{displaymath} 
Now, for $n \geq 2$ and $1/ 16 \pi (n + 1) \leq h < 1 / 16 \pi n$, one gets hence:
\begin{align*} 
\rho_\psi (h) 
& \leq m (\{ z \in \T \, ; \ \Im \psi (z) > 4 \pi n \}) \\
& = \omega_\Omega (a, \{w \in \partial \Omega \, ; \  \Im w > 4 \pi n\}) \leq \omega_\Omega (a, \partial \Omega \setminus \partial \Omega_n) 
\leq \delta_n \leq \delta (h) \,,
\end{align*} 
proving Lemma~\ref{end of the world}.
\qed


\bigskip

\noindent
{\rm Daniel Li}, Univ Lille Nord de France, U-Artois, \\
Laboratoire de Math\'ematiques de Lens EA~2462 \& 
F\'ed\'eration CNRS Nord-Pas-de-Calais FR~2956, \\
Facult\'e des Sciences Jean Perrin, Rue Jean Souvraz, S.P.\kern 1mm 18, \\
F-62\kern 1mm 300 LENS, FRANCE \\ 
daniel.li@euler.univ-artois.fr
\medskip

\noindent
{\rm Herv\'e Queff\'elec}, Univ Lille Nord de France, \\
USTL, Laboratoire Paul Painlev\'e U.M.R. CNRS 8524, \\
F\'ed\'eration CNRS Nord-Pas-de-Calais FR~2956, \\
F-59\kern 1mm 655 VILLENEUVE D'ASCQ Cedex, 
FRANCE \\ 
Herve.Queffelec@univ-lille1.fr
\smallskip

\noindent
{\rm Luis Rodr{\'\i}guez-Piazza}, Universidad de Sevilla, \\
Facultad de Matem\'aticas, Departamento de An\'alisis Matem\'atico \& IMUS,\\ 
Apartado de Correos 1160,\\
41\kern 1mm 080 SEVILLA, SPAIN \\ 
piazza@us.es\par

\end{document}